\documentclass[twoside]{article}
\usepackage{amsmath,amssymb,amscd,amsthm,verbatim,alltt,amsfonts,array}
\usepackage{mathrsfs}
\usepackage[english]{babel}
\usepackage{latexsym}
\usepackage{amssymb}
\usepackage{euscript}
\usepackage{graphicx}

 \def\NZQ{\Bbb}               % the font for N,Z,Q,R,C

 \def\ZZ{{\NZQ Z}}

 \def\C{{\mathcal C}}

\def\Xb{{\bold X}}
 \def\ab{{\bold a}}
 \def\bb{{\bold b}}
 \def\xb{{\bold x}}

 \def\cb{{\bold c}}

 \def\opn#1#2{\def#1{\operatorname{#2}}} % to make operators
 \opn\chara{char} \opn\length{\ell} \opn\pd{pd} \opn\rk{rk}
 \opn\projdim{proj\,dim} \opn\injdim{inj\,dim} \opn\rank{rank}
 \opn\depth{depth} \opn\grade{grade} \opn\height{height}
 \opn\embdim{emb\,dim} \opn\codim{codim}
 
 \opn\Tr{Tr} \opn\bigrank{big\,rank}
 \opn\superheight{superheight}\opn\lcm{lcm}
 \opn\trdeg{tr\,deg}%\emph{
 \opn\reg{reg} \opn\lreg{lreg} \opn\ini{in} \opn\lpd{lpd}
 \opn\size{size} \opn\sdepth{sdepth}
 \opn\link{link}\opn\fdepth{fdepth}\opn\lex{lex}
 \opn\astab{astab}
   \opn\dstab{dstab}
 \opn\div{div} \opn\Div{Div} \opn\cl{cl} \opn\Cl{Cl}
 \opn\Spec{Spec} \opn\Supp{Supp} \opn\supp{supp} \opn\Sing{Sing}
 \opn\Ass{Ass} \opn\Min{Min}\opn\Mon{Mon}
 \opn\astab{astab}
 \opn\Ann{Ann} \opn\Rad{Rad} \opn\Soc{Soc}
 \opn\Im{Im} \opn\Ker{Ker} \opn\Coker{Coker} \opn\Am{Am}
 \opn\Hom{Hom} \opn\Tor{Tor} \opn\Ext{Ext} \opn\End{End}
 \opn\Aut{Aut} \opn\id{id}
 
 \opn\nat{nat}
 \opn\pff{pf}%   \pf exists already
 \opn\Pf{Pf} \opn\GL{GL} \opn\SL{SL} \opn\mod{mod} \opn\ord{ord}
 \opn\Gin{Gin} \opn\Hilb{Hilb}\opn\sort{sort}
 \opn\Tot{Tot}
 \opn\astab{astab}
 \opn\aff{aff} \opn
 \con{conv} \opn\relint{relint} \opn\st{st}
 \opn\lk{lk} \opn\cn{cn} \opn\core{core} \opn\vol{vol}
 \opn\link{link} \opn\star{star}\opn\lex{lex}\opn\set{set}
 \opn\gr{gr}
 
 \def\pot#1#2{#1[\kern-0.28ex[#2]\kern-0.28ex]}

 \opn\dirlim{\underrightarrow{\lim}}
 \opn\inivlim{\underleftarrow{\lim}}
 \let\sect=\cap

 \let\Sect=\bigcap

 \def\Implies{\ifmmode\Longrightarrow \else
         \unskip${}\Longrightarrow{}$\ignorespaces\fi}
 \def\implies{\ifmmode\Rightarrow \else
         \unskip${}\Rightarrow{}$\ignorespaces\fi}
 \def\iff{\ifmmode\Longleftrightarrow \else
         \unskip${}\Longleftrightarrow{}$\ignorespaces\fi}

 \let\:=\colon
 \newtheorem{Theorem}{Theorem}[section]
 \newtheorem{Lemma}[Theorem]{Lemma}
 \newtheorem{Corollary}[Theorem]{Corollary}

 \newtheorem{Example}[Theorem]{Example}
 
 \newtheorem{Definition}[Theorem]{Definition}

 \let\epsilon\varepsilon
 \let\kappa=\varkappa
 \def\qed{\ifhmode\textqed\fi
       \ifmmode\ifinner\quad\qedsymbol\else\dispqed\fi\fi}
 \def\textqed{\unskip\nobreak\penalty50
        \hskip2em\hbox{}\nobreak\hfil\qedsymbol
        \parfillskip=0pt \finalhyphendemerits=0}
 \def\dispqed{\rlap{\qquad\qedsymbol}}
 
 \opn\dis{dis}
 \def\pnt{{\raise0.5mm\hbox{\large\bf.}}}
 
 \opn\Lex{Lex}

\newcommand{\salt}{\vspace*{2.5mm}}     %Maia, 06.05.98

\def\boxit#1{\vbox{\hrule\hbox{\vrule\kern2pt
 \vbox{\kern4pt#1\kern2pt}\kern2pt\vrule}\hrule}}
\def\qed{\hfill
 $\hskip0.3cm
 \boxit{\hsize 2pt \vsize 10pt}$\bigskip\noindent}
\def\C{{\mathchoice {\setbox0=\hbox{$\displaystyle\rm C$}\hbox{\hbox
to0pt{\kern0.4\wd0\vrule height0.9\ht0\hss}\box0}}
{\setbox0=\hbox{$\textstyle\rm C$}\hbox{\hbox
to0pt{\kern0.4\wd0\vrule height0.9\ht0\hss}\box0}}
{\setbox0=\hbox{$\scriptstyle\rm C$}\hbox{\hbox
to0pt{\kern0.4\wd0\vrule height0.9\ht0\hss}\box0}}
{\setbox0=\hbox{$\scriptscriptstyle\rm C$}\hbox{\hbox
to0pt{\kern0.4\wd0\vrule height0.9\ht0\hss}\box0}}}}
\def\Q{{\mathchoice {\setbox0=\hbox{$\displaystyle\rm
Q$}\hbox{\raise
0.15\ht0\hbox to0pt{\kern0.4\wd0\vrule height0.8\ht0\hss}\box0}}
{\setbox0=\hbox{$\textstyle\rm Q$}\hbox{\raise
0.15\ht0\hbox to0pt{\kern0.4\wd0\vrule height0.8\ht0\hss}\box0}}
{\setbox0=\hbox{$\scriptstyle\rm Q$}\hbox{\raise
0.15\ht0\hbox to0pt{\kern0.4\wd0\vrule height0.7\ht0\hss}\box0}}
{\setbox0=\hbox{$\scriptscriptstyle\rm Q$}\hbox{\raise
0.15\ht0\hbox to0pt{\kern0.4\wd0\vrule height0.7\ht0\hss}\box0}}}}
\def\T{{\mathchoice {\setbox0=\hbox{$\displaystyle\rm
T$}\hbox{\hbox to0pt{\kern0.3\wd0\vrule height0.9\ht0\hss}\box0}}
{\setbox0=\hbox{$\textstyle\rm T$}\hbox{\hbox
to0pt{\kern0.3\wd0\vrule height0.9\ht0\hss}\box0}}
{\setbox0=\hbox{$\scriptstyle\rm T$}\hbox{\hbox
to0pt{\kern0.3\wd0\vrule height0.9\ht0\hss}\box0}}
{\setbox0=\hbox{$\scriptscriptstyle\rm T$}\hbox{\hbox
to0pt{\kern0.3\wd0\vrule height0.9\ht0\hss}\box0}}}}
\def\Z{{\mathchoice {\hbox{$\sf\textstyle Z\kern-0.4em Z$}}
{\hbox{$\sf\textstyle Z\kern-0.4em Z$}}
{\hbox{$\sf\scriptstyle Z\kern-0.3em Z$}}
{\hbox{$\sf\scriptscriptstyle Z\kern-0.2em Z$}}}}

\newcommand{\eqnoone}  %normala (fara section) - un numar
   {}
\newcommand{\eqnotwo}  %cu sectiuni - doua numere
   {}

\newcounter{alf}

\newcommand{\adresa}[1]{\par\vspace*{-11pt}
                        \begin{flushright}
                        {\small
                        #1}
                        \end{flushright}
                        }

\begin{document}

  \vskip 1.2 true cm
\setcounter{page}{1}

 \begin{center} {\bf  Generalized mixed product ideals whose powers have a linear resolution} \\

          \medskip

{\sc Monica La Barbiera and Roya Moghimipor}
\end{center}

\begin{abstract}
In this paper we study classes of monomial ideals for which all of its powers have a linear resolution.
Let $K[x_{1},x_{2}]$ be the polynomial ring in two variables over the field $K$, and let $L$ be the generalized mixed product ideal induced by a monomial ideal $I$. It is shown that, if $I\subset K[x_1,x_2]$ and the ideals substituting the monomials in $I$ are Veronese type ideals, then $L^{k}$ has a linear resolution for all $k\geq 1$.
Furthermore, we compute some algebraic invariants of generalized mixed product ideals induced by a transversal polymatroidal ideal.
\end{abstract}

\begin{quotation}
\noindent{\bf Key Words}: {Linear resolution, Monomial localization, Monomial ideals.}

\noindent{\bf 2010 Mathematics Subject Classification}:  13C13, 13D02
\end{quotation}

\thispagestyle{empty}

\section{Introduction}\label{Sec1}

The study of algebraic and homological properties of powers of a monomial ideal has been one of the main topics of combinatorial commutative algebra in
recent years.
Restuccia and Villarreal \cite{RV} introduced mixed product ideals, which form  a particular class of squarefree monomial ideals, and they classified those among these ideals which are normal. In other words, they characterized the mixed product ideals whose Rees ring is normal.

In this paper we consider generalized mixed product ideals which were introduced
by Herzog and Yassemi and which also include the so-called expansions of monomial ideals.
The main objective of this paper is to study all powers of generalized mixed product ideals with linear resolution.

Let $K$ be a field and $K[x_1,\dots,x_n]$ the polynomial
ring in $n$ variables over $K$ with each $x_i$ of degree $1$.
Let $I\subset S$ be a monomial ideal and $G(I)$ its unique
minimal monomial generators.

Let $K$ be a field and  $S = K[x_1, \ldots,  x_n, y_1, \ldots, y_m]$  be the polynomial ring over $K$ in the variables $x_i$ and $y_j$.
A mixed product ideal is a monomial ideal of the form $(I_qJ_r + I_{q'}J_{r'})S$, where for integers $a$ and $b$, the ideal  $I_a$  (resp.\ $J_b$) is the ideal generated by all squarefree monomials of degree $a$ in the polynomial ring $K[x_1,\ldots, x_n]$
(resp.\ of degree $b$ in the polynomial ring $K[y_1, \ldots, y_m]$), and where $0 < q' < q \leq n$, $0 < r < r' \leq m$.
Thus, the ideal $L=(I_qJ_r + I_{q'}J_{r'})S$ is obtained from the monomial ideal $I=(x^qy^r,x^{q'}y^{r'})$ by replacing $x^q$ by $I_q$, $x^{q'}$ by $I_{q'}$, $y^r$ by $J_r$ and $y^{r'}$ by $J_{r'}$.

Together with  Herzog and  Yassemi \cite{HMY} I introduced the generalized mixed product ideals, which  are a far reaching generalization of  the mixed product ideals introduced by Restuccia and Villarreal,  and also generalizes the expansion construction by Bayati and Herzog \cite{BH}.
For this construction we choose for each $i$ a set of new variables $x_{i1},x_{i2},\ldots,x_{im_i}$ and replace each of the factor $x_i^{a_i}$ in each minimal generator $x_1^{a_1}x_2^{a_2}\cdots x_n^{a_n}$ of the monomial ideal $I$ by a monomial ideal in $T_i=K[x_{i1},x_{i2},\ldots,x_{im_i}]$ generated in degree $a_i$.

The class of generalized mixed polymatroidal ideals is a special class of generalized mixed product ideals which for each $i$ we replace each factor $x_i^{a_i}$ in each minimal generator $x_1^{a_1}x_2^{a_2}\cdots x_n^{a_n}$ of $I$ by a polymatroidal ideal in $T_i$ generated in degree $a_i$.
In \cite{M} we showed that powers of a generalized mixed polymatroidal ideal is generalized mixed polymatroidal ideal and monomial localizations of a generalized mixed polymatroidal ideal at monomial prime ideals is again generalized mixed polymatroidal ideal.
A great deal of knowledge on the generalized mixed product ideal is accumulated in several papers \cite{HMY, MM, Mo, MB, MN, MP}.

The present paper is organized as follows. In Section~\ref{two} we study all powers of generalized mixed product ideals with linear resolution.
In Theorem \ref{dim} we compute the dimension of $L$ in the case that $I\subset K[x_1,x_2]$ and the ideals who substitute the generators of $I$
are the ideals $I_{ih_{i}}$ of Veronese type of degree $h_{i}$ in the variables $x_{i1},\dots ,x_{im_{i}}$.
We also characterize the unmixed generalized mixed product ideals, see Theorem \ref{unmixed}.

Our main result (Theorem \ref{power}) says that, if $I\subset K[x_1,x_2]$ is a monomial ideal and the ideals who substitute the generators of $I$ are the ideals of Veronese type generated on degree $h_{i}$ by the set $\{x_{i1}^{a_{i1}}\cdots x_{im_{i}}^{a_{im_{i}}}\mid \sum_{j=1}^{m_{i}}a_{ij}=h_{i},\quad  0\leq a_{ij} \leq 2\}$, then $L^{k}$ has a linear resolution for all $k\geq 1$.

Moreover, we compute the Castelnuovo-Mumford regularity of $L^{k}$ for all $k\geq 1$, see Theorem \ref{reg}.

In Section~\ref{three} we obtain an irredundant primary decomposition for any power of a generalized mixed product ideal $L$ induced by the transversal polymatroidal ideal $I$, see Theorem \ref{decomposition}.
Let $\astab(I)$ and $\dstab(I)$ be the smallest integer $k$
for which $\Ass_{S}(I^{k})$ and $\depth(S/I^{k})$ stabilize, respectively.
In Theorem \ref{assiciated prime of transversal} we compute the stable set of associated prime ideals of $L$, and show that $\astab(L)=1$ if $I$ is a transversal polymatroidal ideal.

Furthermore, we study the depth of powers of $L$, and the index of depth stability of $L$ is determined.
In Theorem \ref{cm} it is shown that $L$ is a Cohen--Macaulay transversal polymatroidal ideal if $I$ is a Cohen-Macaulay transversal polymatroidal ideal.

\section{Generalized mixed product ideals}
\label{two}
Let $R$ be a Noetherian ring and $M$ a finitely generated $R$-module. A prime ideal $P\subset R$ is called an {\em associated prime ideal} of $M$, if
there exists an element $x\in M$ such that $P=\Ann(x)$. Here $\Ann(x)$ is the {\em annihilator} of $x$, that is to say, $\Ann(x):=\{f\in R: fx=0\}$.
The set of associated prime ideals of $M$ is denoted $\Ass_{R}(M)$.

Recall that an ideal $I$ in a Noetherian ring $R$ is $P$-primary, if $\Ass_{R}(R/I)=\{P\}$. In an abuse of notation, one often writes $\Ass_{R}(I)$ instead of $\Ass_{R}(R/I)$.

In \cite{HMY} we introduced the generalized mixed product ideals. Let $S=K[x_{1},\dots,x_{n}]$ be the polynomial ring over a field $K$ in the variables $x_1,\ldots,x_n$ with the maximal ideal $\mathfrak{m}=(x_1,\dots,x_n)$. Let $I\subset S$ be a monomial ideal with $I\neq S$ whose minimal set of generators  is $G(I)=\{\xb^{\ab_1},\ldots, \xb^{\ab_m}\}$. Here we set $\xb^{\ab}=x_1^{\ab(1)}x_2^{\ab(2)}\cdots x_n^{\ab(n)}$ for $\ab=(\ab(1),\ldots,\ab(n))\in\ZZ^{n}_{+}$.
In addition, we consider the polynomial ring $T$ over $K$ in the variables
$x_{11},\ldots,x_{1m_1},x_{21},\ldots,x_{2m_2},\ldots,x_{n1},\ldots,x_{nm_n}$.

For $i=1,\ldots,n$ and $j=1,\ldots,m$ let $L_{i, \ab_j(i)}$ be a monomial ideal in the variables $x_{i1},x_{i2},\ldots,x_{im_i}$ such that
\begin{eqnarray}
\label{inclusion}
L_{i, \ab_j(i)}\subset L_{i, \ab_k(i)} \quad \text{whenever} \quad \ab_j(i)\geq \ab_k(i).
\end{eqnarray}

Given these ideals we define for $j=1,\ldots,m$ the monomial ideals
\begin{eqnarray}
\label{lj}
L_j=\prod_{i=1}^nL_{i, \ab_j(i)} \subset T,
\end{eqnarray}
and set $L=\sum_{j=1}^mL_j$.
The ideal $L$ is called a {\em generalized mixed product ideal} induced by $I$.

\begin{Example}
\label{twovariables}
{\em
Let $T=K[x_{11},\ldots, x_{1m_{1}},x_{21}, \ldots, x_{2m_{2}}]$, and let $L=I_{q}J_{r}+I_{p}J_{s}\subset T$ be the generalized mixed product ideal
induced by a monomial ideal $I=(x^qy^r,x^py^s)$, where for integers $\alpha$ and $\beta$, the ideal  $I_{\alpha}$  (resp.\ $J_{\beta}$) is the ideal generated by all squarefree monomials of degree $\alpha$ in the polynomial ring $K[x_{11},\ldots, x_{1m_{1}}]$
(resp.\ of degree $\beta$ in the polynomial ring $K[x_{21}, \ldots, x_{2m_{2}}]$), and where $0 < p < q \leq m_{1}$, $0 < r < s \leq m_{2}$.
Ideals of this type are called squarefree Veronese ideals.
}
\end{Example}

We want to compute the dimension of $L$ in the case that $I\subset K[x_1,x_2]$ and the ideals who substitute the generators of $I$
are the ideals $I_{ih_{i}}$ of Veronese type of degree $h_{i}$ in the variables $x_{i1},\dots ,x_{im_{i}}$.

Let $\mathcal{F}$ be a non empty subset of $\{x_{11},\ldots,x_{1m_1},x_{21},\ldots,x_{2m_2},\ldots,x_{n1},\ldots,x_{nm_n}\}$. We denote by $P_{\mathcal{F}}$ the prime ideal of $T$ generated by the variables whose index is in $\mathcal{F}$.

Herzog and Hibi \cite{HH} proposed the concept of discrete polymatroid and they studied some algebraic and combinatorics properties related to it.
Let $I$ be a monomial ideal of $S$ which is generated in a single degree.
Then $I$ is said to be {\em polymatroidal}, if for any
two elements $u,v\in G(I)$ such that $\deg_{x_i}u>\deg_{x_i}v$
there exists an index $j$ with $\deg_{x_j}u<\deg_{x_j}v$
such that $x_j(u/x_i)\in I$.
The Veronese type ideals is a special class of polymatroidal ideals, introduced in \cite{V}.

\begin{Definition}
\label{Veronese}
{\em
Given a sequence of integers $(r_{1},\dots,r_{n};d)$ such that $1 \leq r_{j} \leq d \leq \sum_{i=1}^{n}r_{i}$ for all $j$, we define $\mathcal{J}$
as the set of partitions
\[
\mathcal{J}=\{(\ab(1),\ldots,\ab(n))\in \ZZ^{n} \rm {such \; \; that
\; \;} \ab(1)+\dots+\ab(n)=d; \quad 0\leq \ab(i) \leq r_{i} \quad
\forall i\},\] and $H$ as the set of monomials $H=\{\xb^{\ab}\mid
\ab\in \mathcal{J}\}=\{f_{1},\dots,f_{q}\}$. The ideal $I=(H)\subset
S$ is said to be of {\em Veronese type} of degree $d$ with defining
sequence $(r_{1},\dots,r_{n};d)$. If $r_{i}=1$ for all $i$ we call
$I$ the $d$th {\em squarefree Veronese ideal}, and if $r_{i}=d$ for
all $i$ we call $I$ the $d$th {\em Veronese ideal}. }
\end{Definition}

In the following, we give a description of the associated prime ideals of $L$ where the ideals $L_{ij}$ are Veronese type ideals.

\begin{Theorem}
\label{Ass}
Let $T=K[x_{11},\ldots,x_{1m_1},x_{21},\ldots,x_{2m_2}]$ with $m_{1},m_{2}> 1$. Let
\[
L=\sum_{1\leq h_{i} \leq 2m_{i},\sum _{i=1}^2 h_{i}=d}I_{1h_{1}}I_{2h_{2}}
\]
where the ideals $I_{ih_{i}}$ are the ideals of Veronese type
generated on degree $h_{i}$ by the set $\{x_{i1}^{a_{i1}}\cdots x_{im_{i}}^{a_{im_{i}}}\mid \sum_{j=1}^{m_{i}}a_{ij}=h_{i},\quad  0\leq a_{ij} \leq 2\}$.
Then $P_{\mathcal{F}}\in \Ass_{T}(L)$ if and only if $|\mathcal{F}|\leq 2$ for $d=2m_{1}+2m_{2}-1$.
\end{Theorem}

\begin{proof}
The ideal $L$ is a generalized mixed product ideal induced by the ideal $I$ generated by the monomials $x_{1}^{h_{1}}x_{2}^{h_{2}}$
with $\sum_{i=1}^2 h_{i}=d$ and $1\leq h_{i} \leq 2m_{i}$.
$T$. For a subset $\mathcal{F}$ we denote by $P_{\mathcal{F}}$ the prime ideal of $T$ generated by the variables whose index is in $\mathcal{F}$.
We replace the set of variables $\{x_{11},\dots,x_{1m_{1}}\}$ with $\{y_{1},\dots,y_{m_{1}}\}$ and $\{x_{21},\dots,x_{2m_{2}}\}$ with
$\{y_{m_{1}+1},\dots, y_{m_{1}+m_{2}}\}$.
Suppose that $P_{\mathcal{F}}\in \Ass_{T}(L)$. Thus there exists a monomial $f\notin L$ such that $L:f=P_{\mathcal{F}}$.
We prove that we can choose such a monomial $f$ of degree $d-1$ such that $L:f=P_{\mathcal{F}}$.

We assume that $f\notin L$, $L:f=P_{\mathcal{F}}$, $\deg(f)\geq d$ and $f=y_{1}^{a_{1}}\cdots y_{m_{1}+m_{2}}^{a_{m_{1}+m_{2}}}$.
Hence, there exists an integer $r\in \{1,2,\dots,m_{1}+m_{2}\}$ such that $a_{r}> 2$.
Since $L:f=P_{\mathcal{F}}$, it follows that $y_{t}f\in L$ for all $t\in \mathcal{F}$ and $y_{t}f\notin L$ for all $t\notin \mathcal{F}$.
Thus, there exists a monomial $v_{t}\in G(L)$ such that $v_{t}| (y_{t}f)$ for all $t\in \mathcal{F}$.
Since $f\notin L$, we conclude that for all $t\in \mathcal{F}$, the variable $y_{t}$ appears in $v_{t}$ with exponent $a_{t}+1$.
This implies that $a_{t}< 2$ for all $t\in \mathcal{F}$. Therefore, $r\notin \mathcal{F}$.

Now let $g=f/y_{r}$. Since $f\notin L$, it follows that $g\notin L$.
We claim that $L:g=P_{\mathcal{F}}$. In fact, since $g$ divides $f$, it follows that $L:g \subseteq L:f$.
Hence, $L:g \subseteq P_{\mathcal{F}}$.
Now, we want to prove that $P_{\mathcal{F}}\subseteq L:g$.
Let $y_{t}\in P_{\mathcal{F}}$. Since $P_{\mathcal{F}}=L:f$, it follows that $y_{t}f\in L$, so there exists
$v_{t}\in G(L)$ such that $v_{t}\mid y_{t}f$.
Now, since $a_{r}-1\geq 2$, we have that $v_{t} \mid y_{t}f/y_{r}$.
Hence, $y_{t}\in L:f/y_{r}=L:g$.
Therefore, $P_{\mathcal{F}}=L:g$

We proceed as before and obtain after a finite number of these reductions a monomial $f\notin L$
of degree $d-1$ such that $P_{\mathcal{F}}=L:f$.
Then $fy_{t}\in L$ for all $t\in \mathcal{F}$ and $fy_{t}\notin L$ for all $t\notin \mathcal{F}$.
Notice that $a_{t}+1\leq 2$ for all $t\in \mathcal{F}$, and $a_{t}\leq 2$ for all $t\notin \mathcal{F}$.
Thus, $a_{t}=2$ for all $t\notin \mathcal{F}$. Hence, $f=\prod _{t\in \mathcal{F}} y_{t}^{a_{t}} \prod _{t\notin \mathcal{F}}y_{t}^{2}$
with $0\leq a_{t}< 2$ for all $t\in \mathcal{F}$. In particular we see that $\deg(\prod _{t\notin \mathcal{F}} y_{t}^{2})=2(m_{1}+m_{2}-|\mathcal{F}|)=q$. Then
\begin{eqnarray*}
2m_{1}+2m_{2}\geq& \sum _{t\in \mathcal{F}}(a_{t} +1)+q
= \sum _{t\in \mathcal{F}}a_{t}+|\mathcal{F}|+q
= \deg(f)+|\mathcal{F}|
= d-1+|\mathcal{F}|.
\end{eqnarray*}
Hence $|\mathcal{F}| \mid \leq 2$.
Conversely, assume that $|\mathcal{F}|\leq 2$ for $d=2m_{1}+2m_{2}-1$. Then for any monomial $v\in G(L)$ there exists an integer $l\in \mathcal{F}$
such that $y_{l}$ divides $v$. Thus, $L \subset P_{\mathcal{F}}$. By using the condition $2m_{1}+2m_{2}\geq d-1+|\mathcal{F}|$ we have
$|\mathcal{F}|+2m_{1}+2m_{2}-2|\mathcal{F}|\geq d-1$, which together with $2m_{1}+2m_{2}-2|\mathcal{F}|\leq d-1$ implies that
there exists an integer $b_{t}< 2$ for all $t\in \mathcal{F}$ such that $b_{t}|\mathcal{F}|+2m_{1}+2m_{2}-2|\mathcal{F}|=d-1$.
Therefore, the monomial $f=\prod _{t\in \mathcal{F}}y_{t}^{b_{t}} \prod _{t\notin \mathcal{F}} y_{t}^{2}$ has degree $d-1$, and hence
$f\notin L$. Thus $P_{\mathcal{F}}\subseteq L:f$.

Suppose, on the contrary, that $P_{\mathcal{F}}$ is a proper subset of $L:f$. Then there exists a monomial $w$, in the variables $y_{t}$ with $t\notin \mathcal{F}$, of degree at least 1 such that $fw\in L$. Thus, we obtain that there exists a monomial
$v=y_{1}^{a_{1}}\cdots y_{m_{1}+m_{2}}^{a_{m_{1}+m_{2}}}\in G(L)$ such that $v$ divides $fw$.
Observe that $a_{t}\leq b_{t}$ for all $t\in \mathcal{F}$, because $w\in K[y_{t}\mid t\notin \mathcal{F}]$.
This implies that
\begin{eqnarray*}
d=\deg(v)=\sum _{t=1}^{m_{1}+m_{2}}a_{t} \leq \sum _{t\in \mathcal{F}} b_{t}+2m_{1}+2m_{2}-2|\mathcal{F}|=\deg(f)=d-1;
\end{eqnarray*}
a contradiction.
Therefore, $P_{\mathcal{F}}$ is not a proper subset of $L:f$. This shows that $P_{\mathcal{F}}\in \Ass_{T}(L)$.
\end{proof}

\begin{Example}
\label{L7-Ass}
{\em
Let $T=K[x_{11},x_{12},x_{21},x_{22}]$ be a polynomial ring over a field $K$, and let
$L=\sum_{1\leq h_{i} \leq 4,\sum _{i=1}^2 h_{i}=7}I_{1h_{1}}I_{2h_{2}}$
where the ideals $I_{ih_{i}}$ are the ideals of Veronese type
generated on degree $h_{i}$ by the set $\{x_{i1}^{a_{i1}}\cdots x_{im_{i}}^{a_{im_{i}}}\mid \sum_{j=1}^{m_{i}}a_{ij}=h_{i},\quad 0\leq a_{ij} \leq 2\}$.
This implies that
\[
L=(x_{11}^{2}x_{12}x_{21}^{2}x_{22}^{2},
x_{11}x_{12}^{2}x_{21}^{2}x_{22}^{2},
x_{11}^{2}x_{12}^{2}x_{21}^{2}x_{22},
x_{11}^{2}x_{12}^{2}x_{21}x_{22}^{2}).
\]
By using Theorem \ref{Ass} one obtains
\[
\Ass_{T}(L)=\{(x_{11}),(x_{12}),(x_{21}),(x_{22}),(x_{11},x_{12}),(x_{21},x_{22}),(x_{11},x_{21}),
(x_{11},x_{22}),
\]
\[
(x_{12},x_{21}),(x_{12},x_{22})\}.
\]
}
\end{Example}

\begin{Definition}
\label{vertex cover}
{\em
A {\em vertex cover} of $L$ is a subset $\mathcal{C}$ of
\[
\{x_{11},\ldots,x_{1m_1},x_{21},\ldots,x_{2m_2},\ldots,x_{n1},\ldots,x_{nm_n}\}
\]
such that each $v\in G(L)$ is divided by some variables of $\mathcal{C}$.
Such a vertex cover is called {\em minimal}
if no proper subset of $\mathcal{C}$ is vertex cover.
}
\end{Definition}

The minimal cardinality of the vertex covers of $L$ is denoted by $\height(L)$. The reader can find more information in \cite{V}.

\begin{Theorem}
\label{dim}
Let $T=K[x_{11},\ldots,x_{1m_1},x_{21},\ldots,x_{2m_2}]$ with $m_{i}=q> 1$ for $i=1,2$. Let
$
L=\sum_{1\leq h_{i} \leq 2m_{i},\sum _{i=1}^2 h_{i}=d}I_{1h_{1}}I_{2h_{2}}
$
where the ideals $I_{ih_{i}}$ are the ideals of Veronese type
generated on degree $h_{i}$ by the set $\{x_{i1}^{a_{i1}}\cdots x_{im_{i}}^{a_{im_{i}}}\mid \sum_{j=1}^{m_{i}}a_{ij}=h_{i},\quad  0\leq a_{ij} \leq 2\}$.
\begin{itemize}
\item[(\.{i})] If $2\leq d \leq 2m_{1}+2$, then $\dim(T/L)=q$.

\item[(\.{i}\.{i})] If $d=2(m_{1}+m_{2})-1$, then $\dim(T/L)=m_{1}+m_{2}-1$.
\end{itemize}
\end{Theorem}

\begin{proof}
(\.{i}) Let $L=\sum_{r=1}^{m}I_{1h_{1_{r}}}I_{2h_{2_{r}}}$ where the ideals $I_{1h_{1_{r}}}$ in
$K[x_{11},x_{12},\ldots,x_{1m_1}]$ and the ideals $I_{2h_{2_{r}}}$ in
$K[x_{21},x_{22},\ldots,x_{2m_2}]$ are Veronese type ideals of degree $h_{1_{r}}$ and $h_{2_{r}}$, respectively.
Assume further that
\[
1\leq h_{1_{1}}< \dots <h_{1_{m_{1}}}\leq 2m_1 \quad  \text{and}  \quad 2m_2\geq h_{2_{1}}>\dots> h_{2_{m_{2}}}\geq1,
\]
where $h_{1_{r+1}}=h_{1_{r}}+1$, $h_{2_{r-1}}=h_{2_{r}}+1$ for all $r=1,\dots,m-1$.

Let $G(L)$ be the unique minimal set of monomial generators of $L$. A vertex cover of $L$ is a subset $\mathcal{C}$ of $\{x_{11},\ldots,x_{1m_1},x_{21},\ldots,x_{2m_2}\}$ such that each $u\in G(L)$ is divided by some variables of $\mathcal{C}$. Thus, $P$ is a minimal prime ideal of $L$ if and only if $P=(\mathcal{C})$ for some minimal vertex cover $\mathcal{C}$ of $L$.

We may suppose that $h_{1_{r}}+h_{2_{r}}=d$ for all $r=1,\dots,m$.
Then the ideal $L$ is a generalized mixed product ideal induced by the ideal
$$I=(x_1^{h_{1_{1}}}x_2^{h_{2_{1}}},\dots ,x_1^{h_{1_{m_{1}}}}x_2^{h_{2_{m_{2}}}}).$$
Then
\[
I=(x_1^{h_{1_{1}}})\sect\Sect_{r=1}^{m-1}(x_1^{h_{1{_{r+1}}}},x_2^{h_{2{_r}}})\sect (x_2^{h_{2_{m_{2}}}}).
\]
This is the irredundant decomposition of $I$ as an intersection of irreducible ideals. It then follows from \cite[Theorem 1.2]{Mo}
that
\[
L=(I_{1h_{1_{1}}})\cap\sum_{r=1}^{m-1}(I_{1h_{1_{r+1}}},I_{2h_{2_{r}}})\cap(I_{2h_{2_{m_{2}}}}).
\]
Hence, \cite[Corollary 1.7]{Mo} implies that
\[
\height(L)=\min_{r=1,\dots,m-1}\{\height(I_{1h_{1_{1}}}),\height(I_{2h_{2_{m_{2}}}}),\height(I_{1h_{1_{r+1}}})+\height(I_{2h_{2_{r}}}
)\}.
\]
Now we consider $2\leq d \leq  2m_{1}+2$.
Then the minimal vertex covers of $L$ are $\mathcal{C}_{1}=\{x_{11},\ldots,x_{1m_1}\}$, $\mathcal{C}_{2}=\{x_{21},\ldots,x_{2m_2}\}$.
Therefore,
\[
\dim(T/L)=\dim T-\height(L)=2q-\min\{m_{1},m_{2}\}.
\]

(\.{i}\.{i}) The minimal vertex covers of $L$ are $\mathcal{C}_{l}=\{x_{1l}\}$ for $l=1,\dots,m_{1}$ and  $\mathcal{C}'_{p}=\{x_{2p}\}$ for $p=1,\dots,m_{2}$.
Hence, the minimal cardinality of the vertex covers of $L$ is $\height(L)=1$, as desired.
\end{proof}

\begin{Example}
\label{L3}
{\em
Let $T=K[x_{11},x_{12},x_{21},x_{22}]$ be a polynomial ring over a field $K$, and let
$
L=\sum_{1\leq h_{i} \leq 2,\sum _{i=1}^2 h_{i}=3}I_{1h_{1}}I_{2h_{2}}\subset T
$
where the ideals $I_{ih_{i}}$ are the ideals of Veronese type
generated on degree $h_{i}$ by the set $\{x_{i1}^{a_{i1}}\cdots x_{im_{i}}^{a_{im_{i}}}\mid \sum_{j=1}^{m_{i}}a_{ij}=h_{i},\quad  0\leq a_{ij} \leq 2\}$.
The ideal $L$ is a generalized mixed product ideal induced by the ideal $I$ generated by the monomials $x_{1}^{h_{1}}x_{2}^{h_{2}}$
with $h_{1}+h_{2}=3$ and $h_{1},h_{2}\geq 1$. Then $I=(x_{1}x_{2}^{2},x_{1}^{2}x_{2})$.
One easily checks that
\[
I=(x_1)\sect (x_{1}^{2},x_{2}^{2})\sect (x_{2}).
\]
This is the irredundant decomposition of $I$ as an intersection of irreducible ideals.
It then follows that the minimal vertex covers of $L$ are $\mathcal{C}_{1}=\{x_{11},x_{12}\}$, $\mathcal{C}_{2}=\{x_{21},x_{22}\}$.
Therefore, $\height(L)=|\mathcal{C}_{1}|=|\mathcal{C}_{2}|=2$, and hence $\dim(T/L)=2$ by Theorem \ref{dim}.
}
\end{Example}

A monomial ideal is said to be {\em unmixed}
if all its minimal vertex covers have the same cardinality.
We recall the one-to-one correspondence between the minimal
vertex covers of an ideal and its minimal prime ideals.

Now we assume that the generalized mixed product ideals induced by a monomial ideal $I\subset K[x_{1},x_{2}]$, and give a criterion for its unmixedness in the case that the ideals $I_{ih_{i}}$ are the ideals of Veronese type of degree $h_{i}$ in the variables $x_{i1},\dots ,x_{im_{i}}$.

\begin{Theorem}
\label{unmixed}
Let $T=K[x_{11},\ldots,x_{1m_1},x_{21},\ldots,x_{2m_2}]$ with $m_{i}=q> 1$ for $i=1,2$. Let
$
L=\sum_{1\leq h_{i} \leq 2m_{i},\sum _{i=1}^2 h_{i}=d}I_{1h_{1}}I_{2h_{2}}
$
where the ideals $I_{ih_{i}}$ are the ideals of Veronese type
generated on degree $h_{i}$ by the set $\{x_{i1}^{a_{i1}}\cdots x_{im_{i}}^{a_{im_{i}}}\mid \sum_{j=1}^{m_{i}}a_{ij}=h_{i},\quad  0\leq a_{ij} \leq 2\}$.
\begin{itemize}
\item[(\.{i})] If $q\geq 3$ and $d=2(m_{1}+m_{2})-1$, then $L$ is unmixed.

\item[(\.{i}\.{i})] If $q=2$ and $2\leq d \leq 2(m_{1}+m_{2})-1$, then $L$ is unmixed.
\end{itemize}
\end{Theorem}

\begin{proof}
(\.{i}) The ideal $L$ is a generalized mixed product ideal induced by the ideal $I$ generated by the monomials $x_{1}^{h_{1}}x_{2}^{h_{2}}$
with $h_{1}+h_{2}=d$ and $h_{1},h_{2}\geq 1$. One has $\mathcal{C}_{p}=\{x_{1p}\}$, $p=1,\dots,m_{1}$, $\mathcal{C}'_{p'}=\{x_{2p'}\}$, $p'=1,\dots,m_{2}$ are the minimal vertex covers of $L$. Hence, the minimal cardinality of the vertex covers of $L$ is one. Then the assertion follows.

(\.{i}\.{i}) The minimal vertex covers of $L$ are $\mathcal{C}_{1}=\{x_{11},x_{12}\}$, $\mathcal{C}_{2}=\{x_{21},x_{22}\}$. Thus, all the minimal vertex covers have the same cardinality, as desired.
\end{proof}

\begin{Example}
\label{L15}
{\em
Let $T=K[x_{11},x_{12},x_{13},x_{14},x_{21},x_{22},x_{23},x_{24}]$ be a polynomial ring over a field $K$, and let
$
L=\sum_{1\leq h_{i} \leq 14,\sum _{i=1}^2 h_{i}=15}I_{1h_{1}}I_{2h_{2}}
$
where the ideals $I_{ih_{i}}$ are the ideals of Veronese type
generated on degree $h_{i}$ by the set $\{x_{i1}^{a_{i1}}\cdots x_{im_{i}}^{a_{im_{i}}}\mid \sum_{j=1}^{m_{i}}a_{ij}=h_{i},\quad 0\leq a_{ij} \leq 2\}$.
The minimal vertex covers of $L$ are $\mathcal{C}_{l}=\{x_{1l}\}$ for $l=1,\dots,4$ and  $\mathcal{C}'_{p}=\{x_{2p}\}$ for $p=1,\dots,4$.
Therefore, $\height(L)=|\mathcal{C}_{l}|=|\mathcal{C}'_{p}|=1$ for $l=1,\dots,4$ and $p=1,\dots,4$. Hence, by applying Theorem \ref{unmixed} we obtain that $L$ is unmixed.
}
\end{Example}

In the following, we study all powers of generalized mixed product ideals induced by a monomial ideal $I\subset K[x_{1},x_{2}]$, where the ideals substituting the monomials in $I$ are Veronese type ideals.
Note that a generalized mixed product ideal depends not only on $I$ but also on the family $L_{ij}$.
Then for a monomial ideal $I$ with $G(I)=\{\xb^{\ab_1},\ldots, \xb^{\ab_m}\}$ we set
$L(I;\{L_{ij}\})=\sum_{j=1}^{m}\prod _{i=1}^{n}L_{i, \ab_j(i)}$.

\begin{Lemma}
\label{linear}
Let $T=K[x_{11},\ldots,x_{1m_1},x_{21},\ldots,x_{2m_2}]$ be a polynomial ring over a field $K$ with $m_{1},m_{2}> 1$. Let
\[
L=\sum_{1\leq h_{i} \leq 2m_{i},\sum _{i=1}^2 h_{i}=d}I_{1h_{1}}I_{2h_{2}}
\]
where the ideals $I_{ih_{i}}$ are the ideals of Veronese type
generated on degree $h_{i}$ by the set $\{x_{i1}^{a_{i1}}\cdots x_{im_{i}}^{a_{im_{i}}}\mid \sum_{j=1}^{m_{i}}a_{ij}=h_{i},\quad  0\leq a_{ij} \leq 2\}$, and let
\[
L'=\sum_{1\leq h'_{i} \leq 2m_{i},\sum _{i=1}^2 h'_{i}=d'}I_{1h'_{1}}I_{2h'_{2}}
\]
where the ideals $I_{ih'_{i}}$ are the ideals of Veronese type
generated on degree $h'_{i}$ by the set $\{x_{i1}^{a'_{i1}}\cdots x_{im_{i}}^{a'_{im_{i}}}\mid \sum_{j=1}^{m_{i}}a'_{ij}=h'_{i},\quad  0\leq a'_{ij} \leq 2\}$.
Assume further that $2\leq d \leq 2m_{1}+2m_{2}-1$ and $2\leq d' \leq 2m_{1}+2m_{2}-1$.
Then $LL'$ has a linear resolution.
\end{Lemma}

\begin{proof}
Let $L(I;\{L_{ij}\})=\sum_{1\leq h_{i} \leq 2m_{i},\sum _{i=1}^2 h_{i}=d}I_{1h_{1}}I_{2h_{2}}$
be the generalized mixed product ideal, where the ideals $I_{ih_{i}}\subset K[x_{i1},\dots,x_{im_{i}}]$ are the ideals of Veronese type generated on degree $h_{i}$ by the set $\{x_{i1}^{a_{i1}}\cdots x_{im_{i}}^{a_{im_{i}}}\mid \sum_{j=1}^{m_{i}}a_{ij}=h_{i},\quad 0\leq a_{ij} \leq 2\}$.
Let $\ZZ_{+}$ be the set of non-negative integers and $\ZZ^{m_{1}+m_{2}}_{+}$ be the set of the vectors $(\ab_{1},\ab_{2})$ with $\ab_{1}\in \ZZ^{m_{1}}_{+}$ and $\ab_{2}\in \ZZ^{m_{2}}_{+}$, i.e., $(\ab_{1},\ab_{2})=(a_{11},\dots,a_{1m_{1}},a_{21},\dots,a_{2m_{2}})\in \ZZ^{m_{1}+m_{2}}_{+}$ with $a_{1j},a_{2l}\geq 0$.
Furthermore, let
\[
\mathcal{A}=\{(\ab_{1},\ab_{2})\in \ZZ^{m_{1}+m_{2}}_{+}: \sum_{j=1}^{m_{1}}a_{1j}+ \sum_{l=1}^{m_{2}}a_{2l}=d , |\ab_{1}|=h_{1}, |\ab_{2}|=h_{2}, \quad  0\leq a_{1j},a_{2l}\leq 2\}
\]
be the set of the vector exponents of the elements of $G(L(I;\{L_{ij}\}))$.
The modulus of the vector $(\ab_{1},\ab_{2})$ is the number
$|(\ab_{1},\ab_{2})|=|\ab_{1}|+|\ab_{2}|=\sum _{j=1}^{m_{1}}a_{1j}+\sum _{l=1}^{m_{2}}a_{2l}$.
We assume that $(\ab_{1},\ab_{2}),(\bb_{1},\bb_{2})\in \mathcal{A}$ with $a_{1r}> b_{1r}$ for some $r$. Hence, $a_{1r}-1\geq b_{1r}$. Then
\[
|(a_{11},\dots,a_{1r}-1,\dots,a_{1m_{1}},a_{21},\dots,a_{2m_{2}})|=|(b_{11},\dots,\dots,b_{1m_{1}},b_{21},\dots,b_{2m_{2}})|-1<
\]
\[
|(b_{11},\dots,\dots,b_{1m_{1}},b_{21},\dots,b_{2m_{2}})|.
\]
Thus, there exists $r''$ such that
$(a_{11},\dots,a_{1r}-1,\dots,a_{1r''}+1,\dots,a_{1m_{1}},a_{21},\dots,a_{2m_{2}})\in \mathcal{A}$.
Notice that $r\neq r''$. Then for $r\neq r''$, one has $a_{1r''}+1\leq b_{1r''}$. Hence, $a_{1r''}< b_{1r''}$.
Then $L(I;\{L_{ij}\})$ is generated by all the monomials $\Xb_{1}^{\ab_{1}}\Xb_{2}^{\ab_{2}}$ with $(\ab_{1},\ab_{2})\in \mathcal{A}$, where
$\Xb_{1}^{\ab_{1}}\Xb_{2}^{\ab_{2}}$ stands for $x_{11}^{a_{11}}\cdots x_{1m_{1}}^{a_{1m_{1}}}x_{21}^{a_{21}}\cdots x_{2m_{2}}^{a_{2m_{2}}}$.

Suppose that $u=\Xb_{1}^{\ab_{1}}\Xb_{2}^{\ab_{2}}$, $v=\Xb_{1}^{\bb_{1}}\Xb_{2}^{\bb_{2}}\in G(L(I;\{L_{ij}\}))$ with
$a_{1r} > b_{1r}$ or $a_{2t}> b_{2t}$, for some $r''\in \{1,\dots,m_{1}\}$ with $a_{1r''} < b_{1r''}$ and $t''\in \{1,\dots,m_{2}\}$ with $a_{2t''}< b_{2t''}$, we have
$(a_{11},\dots,a_{1r}-1,\dots,a_{1r''}+1,\dots,a_{1m_{1}},a_{21},\dots,a_{2m_{2}})\in \mathcal{A}$ and
$(a_{11},\dots,a_{1m_{1}},a_{21},\dots,a_{2t}-1,\dots,a_{2t''}+1,\dots,a_{2m_{2}})\in \mathcal{A}$.
Therefore, $x_{2r''}(u/x_{1r})\in G(L(I;\{L_{ij}\}))$ or $x_{2t''}(u/x_{2t})\in G(L(I;\{L_{ij}\}))$.

The ideal $L(I;\{L_{ij}\})$ is a generalized mixed product ideal induced by the ideal $I$ generated by the monomials $x_{1}^{h_{1}}x_{2}^{h_{2}}$ with $\sum_{i=1}^2 h_{i}=d$, $1\leq h_{i} \leq 2m_{i}$
and the ideal $L(I';\{L_{ij}\})$ is a generalized mixed product ideal induced by the ideal $I'$ generated by the monomials $x_{1}^{h'_{1}}x_{2}^{h'_{2}}$ with $\sum_{i=1}^2 h'_{i}=d$ and $1\leq h'_{i} \leq 2m_{i}$.
Then
\begin{eqnarray*}\label{contained}
I_{i h_{i}}I_{i h'_{i}}\subset I_{iq_{i}}I_{i q'_{i}} \text{ whenever $h_{i}+h'_{i}\geq q_{i}+q'_{i}$}.
\end{eqnarray*}
Thus, \cite[Lemma 2.2]{M} yields that $L(II';\{L_{ij}\})$ is a generalized mixed product ideal, and $L(II';\{L_{ij}\})=L(I;\{L_{ij}\})L(I';\{L_{ij}\})$.
It then follows that the ideal $L(II';\{L_{ij}\})$ is induced by the ideal $II'$ generated by the monomials $x_{1}^{h_{1}+h'_{1}}x_{2}^{h_{2}+h'_{2}}$ with $1\leq h_{i}\leq 2m_{i}$ and $1\leq h'_{i}\leq 2m_{i}$.
Since the ideals $I$ and $I'$ are polymatroidal ideals, it follows that $II'$ is a polymatroidal ideal, by \cite[Theorem 12.6.3]{HH}.
Therefore, \cite[Corollary 12.6.4]{HH} together with \cite[Proposition 8.2.1]{HH} yield
$II'$ has a linear resolution. Hence by \cite[Theorem 2.5]{M} the desired follows.
\end{proof}

\begin{Definition}
The ideal $L\subset T$ has {\em linear quotients} if there is an ordering $u_1, \ldots, u_r$
of the monomials belonging to $G(L)$ with $\deg(u_1)\leq \cdots \leq \deg(u_r)$ such that for
each $2 \leq j \leq r$, the colon ideal $(u_1,\dots, u_{j-1}) : u_j$ is generated by a subset of
\[
\{x_{11},\ldots,x_{1m_1},x_{21},\ldots,x_{2m_2},\ldots,x_{n1},\ldots,x_{nm_n}\}.
\]
\end{Definition}

The reader can find more information in \cite[Definition 6.3.45]{V}.
We are now come to the main result of present paper.

\begin{Theorem}
\label{power}
Let $T=K[x_{11},\ldots,x_{1m_1},x_{21},\ldots,x_{2m_2}]$ with $m_{1},m_{2}> 1$, and let
\[
L=\sum_{1\leq h_{i} \leq 2m_{i},\sum _{i=1}^2 h_{i}=d}I_{1h_{1}}I_{2h_{2}}
\]
where the ideals $I_{ih_{i}}$ are the ideals of Veronese type
generated on degree $h_{i}$ by the set $\{x_{i1}^{a_{i1}}\cdots x_{im_{i}}^{a_{im_{i}}}\mid \sum_{j=1}^{m_{i}}a_{ij}=h_{i},\quad  0\leq a_{ij} \leq 2\}$.
Furthermore, let $2\leq d\leq 2m_{1}+2m_{2}-1$.
Then $L^{k}$ has a linear resolution for all $k\geq 1$.
\end{Theorem}

\begin{proof}
We proceed with induction on $k$.
For $k=1$, we suppose that $L(I;\{L_{ij}\})=\sum_{1\leq h_{i} \leq 2m_{i},\sum _{i=1}^2 h_{i}=d}I_{1h_{1}}I_{2h_{2}}$
be the generalized mixed product ideal, where the ideals $I_{ih_{i}}\subset K[x_{i1},\dots,x_{im_{i}}]$ are the ideals of Veronese type of degree $h_{i}$.
Assume further that $v\in G(L(I;\{L_{ij}\}))$.

We define $\mathcal{I}=(w\in G(L(I;\{L_{ij}\})): w\prec v)$ with $\prec$ the lexicographical order on $x_{11},\dots,x_{1m_{1}},x_{21},\dots,x_{2m_{2}}$ induced by $x_{11}\succ \cdots \succ x_{1m_{1}} \succ x_{21}\succ \cdots \succ x_{2m_{2}}$. We denote by $(v,w)$ the greatest common divisor of $v$ and $w$. It is well-known that $\mathcal{I}:v=(w/(v,w)\mid w \in \mathcal{I})$.

We claim that for any $w\prec v$, there exists a variable
of $T$ in $\mathcal{I}: v$ such that it divides $w/(v,w)$.  Suppose that $v=\Xb_{1}^{\ab_{1}}\Xb_{2}^{\ab_{2}}$, $w=\Xb_{1}^{\bb_{1}}\Xb_{2}^{\bb_{2}}\in G(L(I;\{L_{ij}\}))$.
Indeed, since $w\prec v$, there exists an integer $l$ such that $a_{1l}> b_{1l}$ and $a_{1q}=b_{1q}$ for $q=1,\dots, l-1$.
Then there exists an integer $l'$ with $a_{1l'}< b_{1l'}$ such that $x_{1l'}(v/x_{1l})\in G(L(I;\{L_{ij}\}))$.

We know that $l < l'$, it follows that $x_{1l'}(v/x_{1l})\in \mathcal{I}$ and hence $x_{1l'}\in \mathcal{I}:v$.
It is then clear that $x_{1l'}$ divides $w/(v,w)$, because the $l'$-th component of the vector exponent of $w/(v,w)$ is given by
\[
b_{1l'}-\min\{b_{1l'},a_{1l'}\}=b_{1l'}-a_{1l'}>0.
\]
If we assume that $a_{1q}=b_{1q}$ for all $q=1,\dots, m_{1}$, $a_{2l}> b_{2l}$ and $a_{2r}=b_{2r}$
for all $r=1,\dots,l-1$, then $x_{2l'}\in \mathcal{I}:v$ and $x_{2l'}$ divides $w/(v,w)$. Therefore, $L(I;\{L_{ij}\})$ has linear quotients, and hence $L(I;\{L_{ij}\})$ has a $d$-linear resolution by \cite[Proposition 8.2.1]{HH}.

For $k=2$ Lemma \ref{linear} implies that $L(I;\{L_{ij}\})^2$ has a $2d$-linear resolution and $L(I^2;\{L_{ij}\})=L(I;\{L_{ij}\})^2$.
Now let $k>2$. Thus, $L(I^k;\{L_{ij}\})=L(I^{k-1}I;\{L_{ij}\})$ and we have as in the case $k=2$ that
$L(I^{k-1}I;\{L_{ij}\})=L(I^{k-1};\{L_{ij}\})L(I;\{L_{ij}\})$. By induction hypothesis and Lemma \ref{linear}, $L(I;\{L_{ij}\})^{k-1}$ has a $d(k-1)$-linear resolution and $L(I;\{L_{ij}\})^{k-1}=L(I^{k-1};\{L_{ij}\})$.
Hence, $L(I;\{L_{ij}\})^{k}$ has a $kd$-linear resolution and
\[
L(I;\{L_{ij}\})^{k}=L(I^{k};\{L_{ij}\}).
\]
Then the desired conclusion follows.
\end{proof}

Together with Theorem \ref{power} the following result is also useful to compute the Castelnuovo-Mumford regularity of powers of $L$ in the case that all $L_{ij}$ are Veronese type ideals.

\begin{Theorem}
\label{reg}
Let $T=K[x_{11},\ldots,x_{1m_1},x_{21},\ldots,x_{2m_2}]$ with $m_{1},m_{2}> 1$, and let
\[
L=\sum_{1\leq h_{i} \leq 2m_{i},\sum _{i=1}^2 h_{i}=d}I_{1h_{1}}I_{2h_{2}}
\]
where the ideals $I_{ih_{i}}$ are the ideals of Veronese type
generated on degree $h_{i}$ by the set $\{x_{i1}^{a_{i1}}\cdots x_{im_{i}}^{a_{im_{i}}}\mid \sum_{j=1}^{m_{i}}a_{ij}=h_{i},\quad  0\leq a_{ij} \leq 2\}$.
Then $\reg(L^{k})=kd$ for all $k\geq 1$.
\end{Theorem}

\begin{Example}
\label{L4}
{\em
Let $T=K[x_{11},x_{12},x_{21},x_{22}]$ be a polynomial ring over a field $K$ in the variables $x_{11},x_{12},x_{21},x_{22}$. Let
$L=\sum_{1\leq h_{i} \leq 4,\sum _{i=1}^2 h_{i}=4}I_{1h_{1}}I_{2h_{2}}$,
where the ideals $I_{ih_{i}}$ are the ideals of Veronese type
generated on degree $h_{i}$ by the set $\{x_{i1}^{a_{i1}}\cdots x_{im_{i}}^{a_{im_{i}}}\mid \sum_{j=1}^{m_{i}}a_{ij}=h_{i},\quad 0\leq a_{ij} \leq 2\}$.
Then
\begin{align*}
L
&=(x_{11}^{2}x_{12}x_{21},x_{11}^{2}x_{12}x_{22}, x_{11}x_{12}^{2}x_{21},x_{11}x_{12}^{2}x_{22}, x_{11}x_{21}^{2}x_{22}, x_{12}x_{21}^{2}x_{22}, x_{11}x_{21}x_{22}^{2},x_{12}x_{21}x_{22}^{2},\\
&x_{11}^{2}x_{21}^{2}, x_{11}^{2}x_{21}x_{22}, x_{11}^{2}x_{22}^{2},x_{12}^{2}x_{21}^{2},x_{12}^{2}x_{22}^{2},x_{12}^{2}x_{21}x_{22},
x_{11}x_{12}x_{21}^{2},x_{11}x_{12}x_{22}^{2},x_{11}x_{12}x_{21}x_{22}).
\end{align*}
Therefore, Theorem \ref{reg} yields $\reg(L^{k})=4k$ for all $k\geq 1$.
}
\end{Example}

\section{Stability of depth and Cohen-Macaulayness of generalized mixed product ideals induced by a transversal polymatroidal ideal}
\label{three}
Let $I\subset S=K[x_1,\ldots,x_n]$ be a monomial ideal as in Section 1 with $G(I)=\{\xb^{\ab_1},\ldots, \xb^{\ab_m}\}$, and $L$ be defined as in $(\ref{lj})$.

The main goal of this section is to study algebraic properties of powers of generalized mixed product ideals induced by a transversal polymatroidal ideal.
We denote the set of monomial prime ideals of $S$ by $\mathcal{P}(S)$.
Let $F$ be a non-empty subset of $[n]$. We denote by $P_{F}$ the monomial prime ideal $(x_{i})_{i\in F}$.
We assume that $L(P_{F};\{L_{ij}\})$ be the generalized mixed product ideal induced by $P_{F}$.
Hence, by applying \cite[Theorem 1.2]{Mo} we obtain that
\[
L(P_{F};\{L_{ij}\})=(L_{i,1})_{i\in F}.
\]
We set $P_{F}'=(L_{i,1})_{i\in F}$, where the ideals $L_{i,1}=(x_{i1},\dots,x_{im_{i}})$ for all $i\in F$.
We denote the set of monomial prime ideals of $T$ by $\mathcal{P}(T)$.

A {\em transversal} polymatroidal ideal is an ideal $I$ of the form
\begin{eqnarray}
\label{transversal}
I=P_{F_{1}}P_{F_{2}}\cdots P_{F_{r}},
\end{eqnarray}
where $F_{1}, \dots, F_{r}$ is a collection of non-empty subsets of $[n]$ with $r\geq1$.
By its definition, the product of transversal polymatroidal ideals is again a transversal polymatroidal ideal.
Then by taking powers of the prime ideal factors of $I$ which appear several times in (3), we have
\begin{eqnarray}
\label{unique}
I=\prod _{t=1}^{s}P_{G_{t}}^{a_{t}} \quad \text{with} \quad  a_{t}\geq1
\end{eqnarray}
where $G_{t}\neq G_{q}$ for $t\neq q$.

In order to determine the associated prime ideals of $L$ induced by the transversal polymatroidal ideal $I$ we will need the following lemma.

\begin{Lemma}
\label{2generalized}
\cite[Lemma 1.2]{M}
Let $L(I;\{L_{ij}\})=\sum_{k=1}^r\prod_{i=1}^{n}L_{i,\ab_k(i)}$
and $L(J;\{L_{ij}\})=\sum_{l=1}^s\prod_{i=1}^{n}L_{i,\bb_l(i)}$ be generalized mixed product ideals, respectively,
induced by the monomial ideals $I$ and
$J$ with $G(I)=\{\xb^{\ab_1},\ldots,\xb^{\ab_r}\}$ and
$G(J)=\{\xb^{\bb_1},\ldots,\xb^{\bb_s}\}$. We assume that
\begin{eqnarray}\label{contained}
L_{i, \ab_k(i)}L_{i, \bb_l(i)}\subset L_{i, \ab_p(i)}L_{i, \bb_q(i)} \text{ whenever $\ab_k+\bb_l\geq \ab_p+\bb_q$}.
\end{eqnarray}
Suppose that $G(IJ)=\{\xb^{\cb_1},\ldots,\xb^{\cb_t}\}$. Then given $\cb_j$, there exist $\ab_k$ and $\bb_l$ such that $\cb_j= \ab_k+\bb_l$.  We set
$L_{i,\cb_j(i)}= L_{i, \ab_k(i)}L_{i, \bb_l(i)}$.
Furthermore, let
$L(IJ;\{L_{ij}\})= \sum_{j=1}^t\prod_{i=1}^nL_{i,\cb_j(i)}.$
Then $L(IJ;\{L_{ij}\})$ is a generalized mixed product ideal, and
\[
L(IJ;\{L_{ij}\})=L(I;\{L_{ij}\})L(J;\{L_{ij}\}).
\]
\end{Lemma}

\begin{Theorem}
\label{transversal}
Let $I=P_{F_{1}}P_{F_{2}}\cdots P_{F_{r}}\subset S$
be a transversal polymatroidal ideal, where
$F_{1}, \dots, F_{r}$ is a collection of non-empty subsets of $[n]$ with $r\geq1$.
Let $L(I;\{L_{ij}\})$ be a generalized mixed product ideal induced by the monomial ideal $I$ satisfying condition (5).
Then $L(I;\{L_{ij}\})$ is a transversal polymatroidal ideal, and
\[
L(I;\{L_{ij}\})=L(P_{F_{1}};\{L_{ij}\})L(P_{F_{2}};\{L_{ij}\})\cdots L(P_{F_{r}};\{L_{ij}\}).
\]
\end{Theorem}

\begin{proof}
Let $I=P_{F_{1}}P_{F_{2}}\cdots P_{F_{r}}$. Thus, Lemma \ref{2generalized} implies that
\begin{eqnarray*}
L(I;\{L_{ij}\})
&=&L(P_{F_{1}}P_{F_{2}}\cdots P_{F_{r}};\{L_{ij}\})\\
&=&L(P_{F_{1}};\{L_{ij}\})L(P_{F_{2}};\{L_{ij}\})\cdots L(P_{F_{r}};\{L_{ij}\}),
\end{eqnarray*}
is a generalized mixed product ideal induced by $I$, where $L(P_{F_{d}};\{L_{ij}\})=(L_{i,1})_{i\in F_{d}}$ for all $d$ with $1\leq d \leq r$.
Since $L_{i,1} \subset K[x_{i1},\ldots,x_{im_i}]$ is a monomial prime ideal, this implies that
$L(I;\{L_{ij}\})\subset T$ is a transversal polymatroidal ideal, as desired.
\end{proof}

By Theorem \ref{transversal}, Lemma \ref{2generalized} together with \cite[Theorem 12.6.2]{HH} yield

\begin{Corollary}
\label{resolution}
Let $I=P_{F_{1}}P_{F_{2}}\cdots P_{F_{r}}\subset S$
be a transversal polymatroidal ideal, where
$F_{1}, \dots, F_{r}$ is a collection of non-empty subsets of $[n]$ with $r\geq1$.
Let $L$ be a generalized mixed product ideal induced by the monomial $I$ satisfying condition (5). Then $L^{k}$ has a linear resolution for all $k\geq 1$.
\end{Corollary}

\begin{Example}
\label{polymatroid}
{\em
Let $F_{1}=\{1,2\}$, $F_{2}=\{1,2,3,4\}$, $F_{3}=\{3,5\}$, $F_{4}=\{4,5\}$ and $I=P_{F_{1}}\cdots P_{F_{4}}$ be the transversal
polymatroidal ideal of $K[x_1, \dots , x_5]$.
Notice that
\[
I=(x_1,x_2)(x_1,x_2,x_3,x_4)(x_3,x_5)(x_4,x_5).
\]
Let
\[
L=(L_{1,1},L_{2,1})(L_{1,1},L_{2,1},L_{3,1},L_{4,1})(L_{3,1},L_{5,1})(L_{4,1},L_{5,1})
\]
be the generalized mixed product ideal induced by $I$ in the polynomial ring
\[
T=K[x_{11},\ldots,x_{1m_1},x_{21},\ldots,x_{2m_2},\ldots,x_{51},\ldots,x_{5m_5}],
\]
where the ideals $L_{i,1}$ are Veronese ideals of degree one in the variables $x_{i1},x_{i2},\ldots,x_{im_i}$
for $i=1,\dots,5$.
Therefore, Corollary \ref{resolution} implies that $L^{k}$ has a linear resolution for all $k\geq 1$.
}
\end{Example}

Now we obtain the index of stability and the stable set of associated primes of $L$ induced by a transversal polymatroidal ideal $I$.
Brodmann showed \cite{B} that there exists an integer $k_1$ such that $\Ass_{S}(I^k)=\Ass_{S}(I^{k_1})$ for all $k\geq k_1$.
The smallest such number is called {\em index of stability} of $I$.
We denote this number by $\astab(I)$. The stable set $\Ass_{S}(I^k)$ is denoted by $\Ass^\infty (I)$.

\begin{Theorem}
\label{assiciated prime of transversal}
Let $L$ be the generalized mixed product ideal induced by a transversal polymatroidal ideal $I\subset S$ satisfying condition (5). Then $\astab(L)=1$, that is
\[
\Ass_{T}(L)=\Ass^\infty (L).
\]
\end{Theorem}

\begin{proof}
Let $I=P_{F_{1}}P_{F_{2}}\cdots P_{F_{r}}\subset S$ be a transversal polymatroidal ideal, where
$F_{1}, \dots, F_{r}$ is a collection of non-empty subsets of $[n]$ with $r\geq1$. We assume that $k\geq1$ be an integer. Furthermore, let $L(I;\{L_{ij}\})$ be the generalized mixed product ideal induced by $I$ satisfying condition (5).

We denote $I(P)$ the monomial localization of a monomial ideal $I$, and by $S(P)$ the polynomial ring
over $K$ in the variables which belong to $P$.
The monomial localization $I(P)$ is the monomial ideal which is obtained from $I$ as the image of the $K$-algebra homomorphism $\varphi:S\rightarrow S(P)$
which $\varphi(x_i)=x_i$ if $x_i\in P$, and $\varphi(x_i)=1$, otherwise.
According to Theorem \ref{transversal} we have
\[
L(I;\{L_{ij}\})=L(P_{F_{1}};\{L_{ij}\})L(P_{F_{2}};\{L_{ij}\})\cdots L(P_{F_{r}};\{L_{ij}\}),
\]
where $L(P_{F_{d}};\{L_{ij}\})=(L_{i,1})_{i\in F_{d}}$ for all $d=1,\dots,r$.
Theorem \ref{transversal} together with Lemma \ref{2generalized} now yields $L(I;\{L_{ij}\})^{k}$ is a transversal polymatroidal ideal induced by $I^{k}$, and $L(I;\{L_{ij}\})^{k}=L(I^{k};\{L_{ij}\})$.
By \cite[Lemma 2.3]{HRV}, we obtain that $P\in \Ass_{T}(L(I;\{L_{ij}\})^{k})$ if and only if $\depth (T(P)/L(I;\{L_{ij}\})^{k}(P))=0$. We know that the localization of a transversal polymatroidal ideal is again a transversal polymatroidal ideal.
Then $L(I;\{L_{ij}\})(P)$ is again a transversal polymatroidal ideal.
Since powers of transversal polymatroidal ideals are again transversal polymatroidal ideal, and since by \cite[Theorem 12.6.2]{HH}
transversal polymatroidal ideals have linear resolutions. Thus, all powers of $L(I;\{L_{ij}\})(P)$ have a linear resolution.
Therefore, \cite[Proposition 2.2]{HRV} yields $\depth (T(P)/L(I;\{L_{ij}\})^{q}(P))=0$ for all $q\geq k $
and hence $P \in \Ass_{T}(L(I;\{L_{ij}\})^{q})$ for all $q\geq k$. Thus, $\Ass_{T}(L(I;\{L_{ij}\})) \subset \Ass^\infty (L(I;\{L_{ij}\}))$.

For the converse inclusion, let $P\in \Ass^\infty (L(I;\{L_{ij}\}))$. Then there exists $k\geq 1$ such that $P\in \Ass_{T}(L(I;\{L_{ij}\})^{k})$.
Hence, by applying \cite[Lemma 2.3]{HRV} and \cite[Corollary 4.5]{HRV} we obtain that $P\in \Ass_{T}(L(I;\{L_{ij}\}))$. This yields the
desired conclusion.
\end{proof}

Next we study the irredundant primary decomposition of all powers of generalized mixed product ideals induced by a transversal polymatroidal ideal.
Let $I=P_{F_{1}}P_{F_{2}}\cdots P_{F_{r}}\subset S$ be a transversal polymatroidal ideal.
Without loss of generality we may assume that $\bigcup_{d=1}^{r}F_{d}=[n]$.

Let $G_{I}$ be a graph attached to the ideal $I$ on the vertex set $\{1, \dots,r\}$ with $\{i,j\}$ an edge of $G_{I}$ if $F_{i} \cap F_{j}\neq \emptyset$.
In \cite[Theorem 4.7]{HRV} it is shown that the associated prime ideals of $I$ corresponded to the trees of the graph $G_{I}$
in the following way: $P_{F}\in \Ass_{S}(I)$ if and only if there exists a tree $\mathcal{T}$ such that $P_{F}=P_{\mathcal{T}}$.

For each subgraph $\mathcal{H}$ of $G_{I}$ we associate the prime ideal $P_{\mathcal{H}}=\sum _{d\in \mathcal{V}(\mathcal{H})} P_{F_{d}}$.
Furthermore, for a monomial prime ideal $P_{\mathcal{H}}$ we set
\[
P_{\mathcal{H}}'=\sum _{d\in \mathcal{V}(\mathcal{H})}P_{F_{d}}'.
\]
Notice that $P_{\mathcal{H}}'$ is a generalized mixed product ideal induced by the monomial prime $P_{\mathcal{H}}$,
where the ideals $L_{i,1}=(x_{i1},\dots,x_{im_{i}})$ for all $i\in F_{d}$.

\begin{Theorem}
\label{decomposition}
Let $I$ be a transversal polymatroidal ideal with $\Ass_{S}(I)=\{P_{1},\dots,P_{q}\}$, and let $L$ be the generalized mixed product ideal induced by $I$, where all $L_{i,1}=(x_{i1},x_{i2},\ldots,x_{im_i})$ for $d=1,\dots, r$.
Consider $\mathcal{T}_{1},\dots, \mathcal{T}_{q}$ maximal trees of $G_{I}$ such that $P_{h}=P_{\mathcal{T}_{h}}$ for all $h=1,\dots,q$.
Then
\[
L^{k}=\bigcap_{h=1}^{q}(P'_{h})^{ka_{h}},
\]
is an irredundant primary decomposition of $L^{k}$ for any $k\geq 1$, where $a_{h}=|\mathcal{V}(\mathcal{T}_{h})|$ for all $h$.
\end{Theorem}

\begin{proof}
For a transversal polymatroidal ideal $I=P_{F_{1}}P_{F_{2}}\cdots P_{F_{r}}\subset S$, let $L(I;\{L_{ij}\})$ be the generalized mixed product ideal induced by $I$.
By Theorem ~\ref{transversal} together with Lemma \ref{2generalized} it follows that  $L(I;\{L_{ij}\})^{k}$ is a transversal polymatroidal ideal induced by $I^{k}$, and $L(I;\{L_{ij}\})^{k}=L(I^{k};\{L_{ij}\})$ for all $k\geq 1$.

Suppose that $\Ass_{S}(I)=\{P_{1}, \dots, P_{q}\}$. Consider $\mathcal{T}_{1},\dots, \mathcal{T}_{q}$ maximal trees of $G_{I}$
such that $P_{h}=P_{\mathcal{T}_{h}}$ for all $h=1,\dots, q$.
Therefore, \cite[Corollary 4.10]{HRV} yields
\[
I^{k}=\bigcap_{h=1}^{q}P_{h}^{ka_{h}},
\]
is an irredundant primary decomposition of $I^{k}$, where $a_{h}=|\mathcal{V}(\mathcal{T}_{h})|$ for all $h$.
Hence, by \cite[Proposition 1.2]{BH} we have
\begin{eqnarray*}
L(I;\{L_{ij}\})^{k}=L(I^{k};\{L_{ij}\})=L(\bigcap_{h=1}^{q}P_{h}^{ka_{h}};\{L_{ij}\})=\bigcap_{h=1}^{q} L(P_{h};\{L_{ij}\})^{ka_{h}},
\end{eqnarray*}
as desired.
\end{proof}

Let $I$ be a monomial ideal of $S$. We assume that $t$ be a variable over $S$. The graded subalgebra
$\mathcal{R}(I):=\bigoplus _{k=0}^{\infty}I^{k}t^{k}$
of $S[t]$ is called the {\em Rees algebra} of $I$.
The $K$-algebra $\mathcal{R}(I)/\mathfrak{m}\mathcal{R}(I)$ is called the {\em fiber ring} and its Krull dimension
is called the {\em analytic spread} of $I$, denoted by $\ell (I)$.

By a theorem of Brodmann \cite{Bro}, $\depth (S/I^{k})$ is constant for large $k$.
We call this constant value the {\em limit depth} of $I$, and denote it by $\lim_{k\rightarrow \infty} \depth (S/I^{k})$.
The smallest integer $k_{0}$ such that $\depth (S/I^{k})=\depth (S/I^{k_{0}})$ for all $k\geq k_{0}$ is called the {\em index
of depth stability} and is denoted by $\dstab(I)$.

In the following, we compute the analytic spread of $L$ induced by a transversal polymatroidal ideal $I$.

\begin{Theorem}
\label{analytic}
Let $L$ be a generalized mixed product ideal induced by the transversal polymatroidal ideal $I$ satisfying condition (5).
Then $\depth (T/L)=\depth (T/L^{k})$ for all $k\geq 1$. In particular, we have
\[
\ell(L)=m_1+\dots+m_m-\depth (T/L).
\]
\end{Theorem}

\begin{proof}
Let $I=P_{F_{1}}P_{F_{2}}\cdots P_{F_{r}}\subset S$ be a transversal polymatroidal ideal. Let $L(I;\{L_{ij}\})$ be a generalized mixed product ideal induced by $I$.
Theorem ~\ref{transversal} with Lemma \ref{2generalized} guarantees that
$L(I;\{L_{ij}\})^{k}$ is a transversal polymatroidal ideal induced by $I^{k}$, and $L(I;\{L_{ij}\})^{k}=L(I^{k};\{L_{ij}\})$ for all $k\geq 1$.
We denote the number of connected components of the graph $G_{I}$ by $c(G_{I})$.
We know that $c(G_{I})=c(G_{I^{k}})$ for any $k\geq1$. Then \cite[Theorem 4.12]{HRV} implies that
$\depth (T/L(I;\{L_{ij}\}))=\depth (T/L(I;\{L_{ij}\})^{k})$.
Therefore,
\[
\lim_{k\rightarrow \infty}\depth (T/L(I;\{L_{ij}\})^{k})=\depth (T/L(I;\{L_{ij}\}))
\]
for all $k\geq1$, and hence
\[
\lim_{k\rightarrow \infty}\depth (T/L(I;\{L_{ij}\})^{k})=m_1+\dots+m_m-\ell(L(I;\{L_{ij}\}))
\]
by \cite[Corollary 3.5]{HRV}. Thus, the desired conclusion follows.
\end{proof}

Let $k$ be a nonnegative integer. We define the $k$-th {\em bracket power}
of $I$, to be the ideal $I^{[k]}$, generated by all monomials $u^{k}$, where
$u\in I$ is a monomial.
In particular, $I^{[0]}=S$, and $I^{[1]}=I$.

Note that if $G(I)=\{\xb^{\ab_1},\ldots, \xb^{\ab_m}\}$, then
$G(I^{[k]})=\{\xb^{k\ab_1},\ldots, \xb^{k\ab_m}\}$.
Notice that $I^{[k]}\subset I^{k}$ for all $k$.
When $k\geq2$, the equality holds, if and only if $I$ is principal. One can easily compute that $(I\cap J)^{[k]}=I^{[k]}\cap J^{[k]}$ for
any monomial ideals $I,J \subset S$.

Let $I=P_{F_{1}}P_{F_{2}}\cdots P_{F_{r}}\subset S$ be a transversal polymatroidal ideal. Let $G_{1},\dots, G_{t}$ be the connected components of $G_{I}$, with $t\geq2$.
Thus, $I=I_{1}\cdots I_{t}$, where $I_{a}=\prod_{d\in \mathcal{V}(G_{a})}P_{F_{d}}$.
Furthermore, $G_{a}=G_{I_{a}}$ for all $a$. Note that
$I=I_{1}\cdots I_{t}=I_{1}\cap \cdots \cap I_{t}$, we know that the ideals
$I_{a}$ are generated in pairwise disjoint sets of variables. Therefore, $I^{[k]}=I_{1}^{[k]} \sect \cdots \sect I_{t}^{[k]}$. Hence, \cite[Lemma 2.3]{MB} implies that $\Ass_{S}(I^{[k]})=\Ass_{S}(I_{1})\cup \cdots \cup \Ass_{S}(I_{t})$, where $I_{a}$ is a transversal polymatroidal ideal with the associated connected graph $G_{I_{a}}$.

In the following, we want to determine the associated prime ideals of the $k$-th bracket power $L^{[k]}$ of generalized mixed product ideals $L$ induced by a transversal polymatroidal ideal $I$.

\begin{Theorem}
\label{bracket of transversal}
Let $I=P_{F_{1}}P_{F_{2}}\cdots P_{F_{r}}$ be a transversal polymatroidal ideal with the set of associated prime ideals $\Ass_{S}(I)=\{P_1,\dots , P_q\}$.
Consider $\mathcal{T}_{1},\dots, \mathcal{T}_{q}$ maximal trees of $G_{I}$ such that $P_{h}=P_{\mathcal{T}_{h}}$ for all $h=1,\dots,q$.
Let $L$ be a generalized mixed product ideal induced by $I$, where all $L_{i,1}=(x_{i1},x_{i2},\ldots,x_{im_i})$ for $d=1,\dots, r$.
Then
\[
\Ass_{T}(L^{[k]})=\{P'_{h}\mid h=1,\dots,q\}.
\]
for all $k\geq 1$.
\end{Theorem}

\begin{proof}
The assertion follows by Theorem \ref{transversal} and \cite[Theorem 2.4]{MB}.
\end{proof}

Now we assume that the generalized mixed product ideals
induced by a transversal polymatroidal ideal, and give a criterion for its Cohen-Macaulayness.

\begin{Theorem}
\label{cm}
Let $L$ be the generalized mixed product ideal induced by a Cohen-Macaulay transversal polymatroidal ideal $I=P_{F_{1}}P_{F_{2}}\cdots P_{F_{r}}$ such that $\height(I)=t$.
Assume further that
\[
L(P_{F_{d}};\{L_{ij}\})=(L_{i,1})_{i\in F_{d}}
\]
where $L_{i,1}=(x_{i1},\dots,x_{im_{i}})$ for all $d$ with $d=1,\dots, r$. If $m_{1}=\dots=m_{n}$ and $\projdim(T/L)=m_{1}+\dots+m_{t}$, then $L$ is a Cohen-Macaulay transversal polymatroidal ideal.
\end{Theorem}

\begin{proof}
Let $I=P_{F_{1}}P_{F_{2}}\cdots P_{F_{r}}\subset S$ be a Cohen-Macaulay transversal polymatroidal ideal.
Then \cite[Theorem 4.7]{HRV} implies that
\[
\dim S/I=n-\min \{ |F_{d}|: d=1,\dots,r\}.
\]
We denote by $t$ the minimal cardinality of a set $F_{d}$, where $d=1,\dots,r$.
Thus, $\dim S/I=n-t$.
Let $L(I;\{L_{ij}\})$ be a generalized mixed product ideal induced by the monomial ideal $I$.
We may assume that $\cup_{d=1}^{r}F_{d}=[n]$, and then it follows from \cite[Theorem 4.12]{HRV} that
$n-t=k-1$, where $k$ represents the number of connected components of $G_{I}$.
We know that $n\geq kt$ it follows that $k-1=n-t\geq kt-t=(k-1)t$.
Therefore we obtain that this inequality is valid either if $k=1$ or $t=1$.
If $k=1$, then $t=n$. This implies that $F_{d}=[n]$ for all $d$ and thus $I=\mathfrak{m}^{r}$, the Veronese ideal.
Otherwise $t=1$ and hence $k=n$. In this case $G_{I}$ has $n$ connected components.
Then $r=n$ and $|F_{d}|=1$ for all $d=1,\dots, r$. This yields that $I$ is a principal ideal.
Thus, Theorem ~\ref{transversal}, \cite[Theorem 12.6.7]{HH} together with \cite[Theorem 2.3]{Mo} yield
$L$ is a Cohen-Macaulay transversal polymatroidal ideal, as desired.
\end{proof}

\begin{Example}
\label{maximal}
{\em
Let $I \subset S$ be Veronese ideal of degree $r$ which is generated by all the monomials in the variables $x_{1},\ldots,x_{n}$ of degree $r$.
Furthermore, let $L(\mathfrak{m}^{r};\{L_{ij}\})$ be the generalized mixed product ideal induced by $I=\mathfrak{m}^{r}$, where the ideals substituting the monomials in $I$ are Veronese ideals. It follows from \cite[Lemma 2.3]{HRV} that $\depth (S/I)=0$.
We know that powers of Veronese ideals are again polymatroidal, see \cite[Theorem 12.6.3]{HH}, and since by \cite[Theorem 12.6.2]{HH} Veronese ideals have linear resolutions, we obtain that all powers of $I$ have a linear resolution.
By \cite[Proposition 1.3]{HRV} we have $\depth (S/I^{k})=0$ for all $k\geq 1$. Thus, Lemma \ref{2generalized} together with \cite[Proposition 1.2]{BH} now yields
$L(\mathfrak{m};\{L_{ij}\})\in \Ass_{T}(L(\mathfrak{m}^{r};\{L_{ij}\})^{k})$ for all $k\geq 1$. Then $(x_{11},x_{12}, \dots,x_{nm_{m}})\in \Ass^{\infty}(L)$.
Theorem \ref{assiciated prime of transversal} together with Theorem \ref{analytic} implies that
\[
\astab(L(\mathfrak{m}^{r};\{L_{ij}\}))=\dstab(L(\mathfrak{m}^{r};\{L_{ij}\}))=1.
\]
}
\end{Example}

\medskip

\salt

\adresa{Department of Mathematics and Informatics,
University of Messina,\\
Viale Ferdinando
Stagno d’Alcontres 31, 98166 Messina, Italy \\
E-mail:{\tt monicalb@unime.it}
}

\salt

\adresa{ Department of Mathematics, Safadasht Branch,\\
Islamic Azad University, Tehran, Iran \\
E-mail: {\tt roya\_moghimipour@yahoo.com}
}

\end{document}